\pdfoutput=1
\documentclass{preprint}

\usepackage{graphicx}
\usepackage{amsmath}
\usepackage{amsfonts}
\usepackage{amssymb}

\title{Boundary Element Analysis with trimmed NURBS and a generalized IGA approach}

\author{Gernot Beer$^{1}$, Benjamin Marussig$^{2}$, J\"urgen Zechner $^{3}$, Christian D\"unser$^{4}$ and Thomas-Peter Fries$^{5}$}

\heading{G. Beer, B. Marussig, J. Zechner, Ch. Duenser and T-P Fries}

\address{$^{1}$ Emeritus professor , TU Graz, Lessingstrasse 25, Graz, Austria \\ and conjoint professor Centre for Geotechnical and Materials Modeling, University of Newcastle, Callaghan, Australia, gernot.beer@tugraz.at\\
$^{2}$ PhD Student, TU Graz , Lessingstrasse 25, Graz, Austria marussig@tugraz.at  \\
$^{3}$ Postdoc, TU Graz , Lessingstrasse 25, Graz, Austria, juergen.zechner@tugraz.at  \\
$^{4}$ Staff scientist,TU Graz Lessingstrasse 25, Graz, Austria, duenser@tugraz.at  \\
$^{5}$ Full Professor, TU Graz Lessingstrasse 25, Graz, Austria, fries@tugraz.at }

\keywords{Boundary Element Method, Isogeometric method.}

\abstract{A novel approach to the simulation with the boundary element method using trimmed NURBS patches is presented. The advantage of this approach is its efficiency and easy implementation. The analysis with trimmed NURBS is achieved by double mapping. The variation of the unknowns on the boundary is specified in a local coordinate system and is completely independent of the description of the geometry. The method is tested on a branched tunnel and the results compared with those obtained from a conventional analysis. The conclusion is that the proposed approach is superior in terms of number of unknowns and effort required. }

\begin{document}

\section{INTRODUCTION}

The boundary element method (BEM) has offered an alternative to the finite element method and has been attractive for certain types of problems, such as those involving an infinite or semi-inifinite domain \cite{BeerSmithDuenser2008b}.

The isogeometric approach \cite{Hughes2005a} has led to renewed interest in the method since it only requires a surface discretization and a direct link can be established with geometric modeling technology, without the need to generate a mesh. 
Using NURBS instead of the traditional Serendipity or Lagrange functions  for describing the variation of boundary values, additional benefits are gained because of their higher continuity and efficient refinement strategies \cite{Beer2013}.

Trimmed NURBS patches have been successfully applied to problems where two solids intersect.  Such models can be created quickly in a CAD program and data exported in IGES format. The exported IGES data contain the description of the boundary of objects via NURBS patches and trimming curves defined in the local coordinate system of a patch. This information can then be used to trim the NURBS patch, i.e. to remove part of the patch surface.

In this paper we present a simple but effective approach to the analysis of trimmed surfaces. To the best of our knowledge the proposed method of double mapping has not been previously published. 
Furthermore we propose to generalize the isogeometric concept by approximating the boundary values (tractions, displacements) with functions that are different from the ones used to describe the geometry. Our motivation for this comes from the fact that the data obtained from CAD programs describing the boundary, may not be suitable for describing the boundary values.

\section{GEOMETRY DEFINITION WITH TRIMMED NURBS PATCHES}
\subsection{NURBS patches}

In CAD programs the geometry is described by NURBS patches which are mapped from a unit square with coordinates $u,v$ to the global coordinates \textbf{x} (x,y,z) by
\begin{equation}
\mathbf{x}=  \sum_{b=1}^{B}\sum_{a=1} ^{A}R_{a,b}^{p,q} (u,v)\mathbf{x}_{a,b}
\label{Gmap}
 \end{equation}
where $\mathbf{x}_{a,b}$ are the coordinates of the control points, $p$ and $q$ are the function orders in $u$ and $v$ direction, $A$ and $B$ are the number of control points in $u$ and $v$ direction and $R_{a,b}^{p,q}$ are tensor products of NURBS functions:
\begin{equation}
R_{a,b}^{p,q}=\frac{N_{a,p}(u)N_{b,q}(v)w_{a,b}}{\sum_{\bar{b}=1}^{B}\sum_{\bar{a}=1}^{A}N_{\bar{a},p}(u)N_{\bar{b},q}(v)w_{\bar{a},\bar{b}}}
 \end{equation}
$N_{a,p}(u)$ and $N_{b,q}(v)$ are B-spline functions of local
coordinates $u$ or $v$ of order $p$ or $q$ ($0$ constant , $1$ linear
, $2$ quadratic etc.) and $w_{a,b}$ are weights. The B-spline
functions are defined by a Knot vector with non-decreasing values of
the local coordinate and a recursive computation which starts at order
0. For example for $N_{a,p}(u)$ with a Knot vector $\Xi_{u}=\left(\begin{array}{cccc}u_{1} & u_{2} & \cdots & u_{A+p+1}\end{array}\right)$ we have for $p=0$
\begin{eqnarray}
N_{a,p}(u)&=&1  \quad  for  \quad  u_{a}\leqslant u < u_{a+1}\\
N_{a,p}(u)&=&0 \quad otherwise
\end{eqnarray}
and for $p$=1,2,3$\cdots$
\begin{equation}
\label{ }
N_{a,p}(u)=\frac{u-u_{p}}{u_{a+p}-u_{a}} \cdot N_{a,p-1} + \frac{u_{a+p+1}-u}{u_{a+p+1}-u_{a+1}} \cdot N_{a+1,p-1} \end{equation}
For a more detailed description of NURBS the reader is referred to \cite{Piegl1997b}.

\begin{figure}[t]
\begin{center}
\includegraphics[scale=0.8]{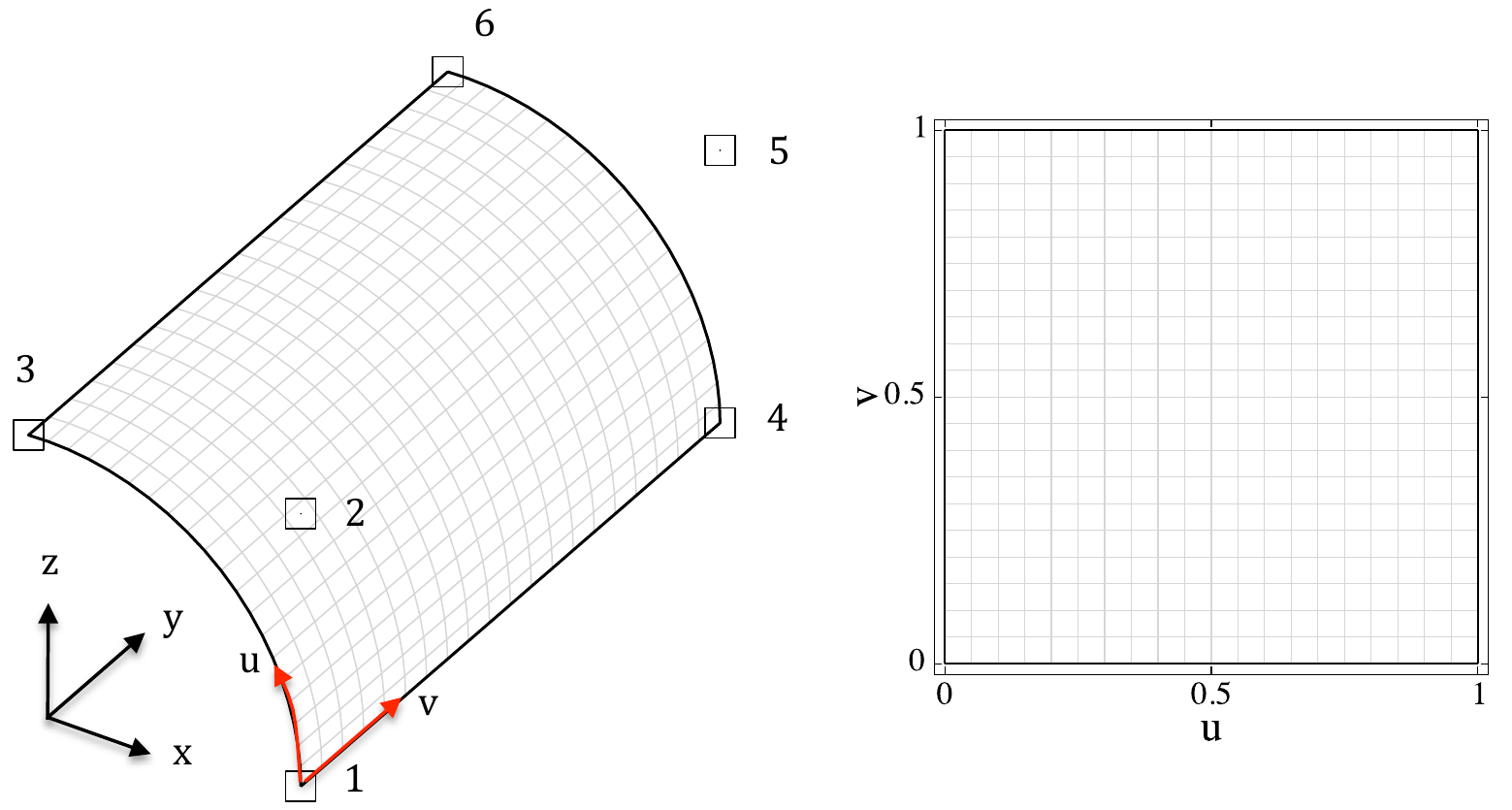}
\caption{Quarter cylinder with control points and local u,v coordinate system}
\label{Cyl}
\end{center}
\end{figure}

As an example we show in Figure \ref{Cyl} the geometrical representation of a quarter cylinder.
For this example the knot vectors in u,v direction are 
\begin{eqnarray}
\Xi_{u} & = & \left(\begin{array}{cccccc}0 & 0 & 0 & 1 & 1 & 1\end{array}\right) \\
\Xi_{v} & = & \left(\begin{array}{cccccc}0 & 0  & 1 & 1\end{array}\right) 
\end{eqnarray}
and the weights are given by
\begin{equation}
\label{ }
\textbf{w}=\left(\begin{array}{ccc}1 & 0.7 & 1 \\1 & 0.7 & 1\end{array}\right)
\end{equation}

\subsection{Analysis with trimmed NURBS patches}

If there is an intersection of NURBS patches the CAD program provides trimming information.
The trimming information comprises one or more trimming curves, which are B-splines and are defined in the local coordinate of the NURBS patch to be trimmed. A method for performing an analysis on trimmed surfaces has already been presented in \cite{Schmidt2012} and involves finding the intersection of the trimming curve with the underlying NURBS surface and a reconstruction of the knot spans and control points.

\begin{figure}
\begin{center}
\includegraphics[scale=1.0]{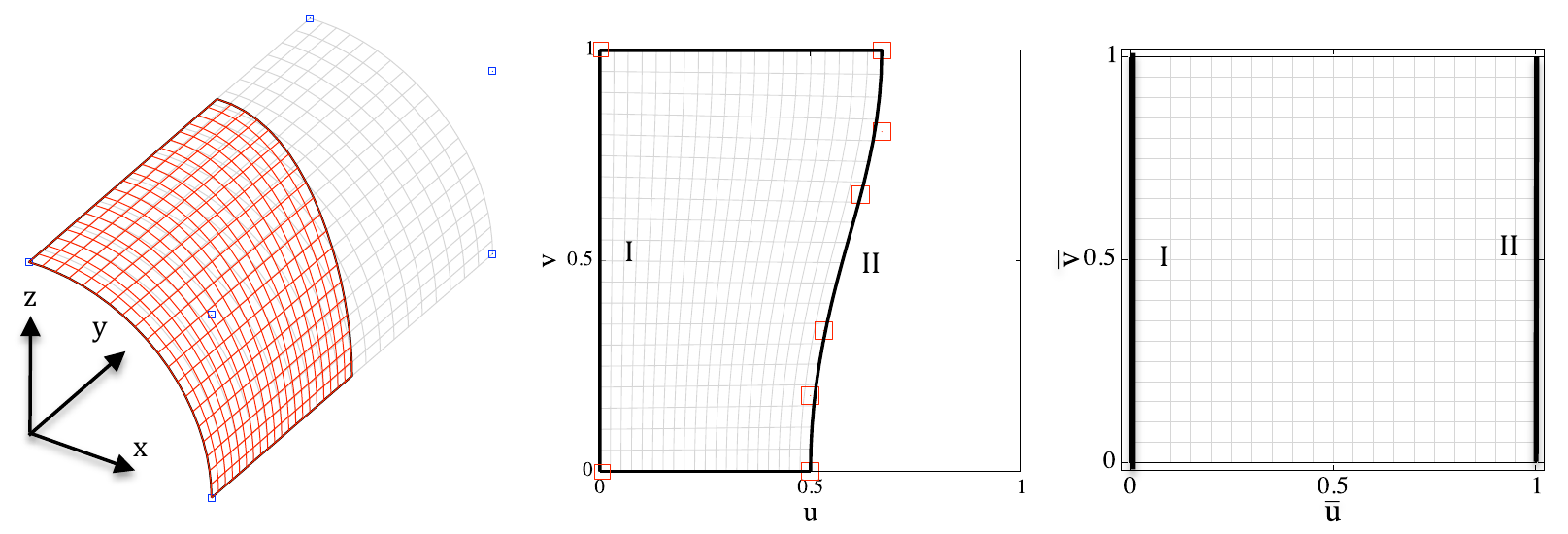}
\caption{Trimming of a quarter cylinder using the double mapping algorithm}
\label{TrimC}
\end{center}
\end{figure}

Here we present a novel approach that appears to be simpler to implement and more efficient, since all it involves is a mapping.  
For explaining the proposed trimming method we use a simple example and assume that the cylinder in Figure \ref{Cyl} is trimmed by 2 trimming curves obtained from the CAD program, marked I and II in Figure \ref{TrimC}. Trimming curve I is a B-spline of order $p=1$ and has 2 control points and  trimming curve II is of order $p=3$ and has 6 control points. The idea is to map the trimmed area from the $u,v$ coordinate system to an $\bar{u}$ ,$\bar{v}$ coordinate system that represents a unit square. The trimming curves map as straight lines along $\bar{v}$ at $\bar{u}=0$ and $\bar{u}=1$ in this system.
The mapping from the $\bar{u}$ ,$\bar{v}$ to the $u,v$ system is given by
\begin{eqnarray}
u & = & N_{1}(\bar{u})u_{I} (\bar{v}) + N_{2}(\bar{u})u_{II} (\bar{v})\\
v & = & N_{1}(\bar{u})v_{I} (\bar{v}) + N_{2}(\bar{u})v_{II} (\bar{v}) 
\end{eqnarray}
where
\begin{eqnarray}
N_{1}(\bar{u})& = &1-\bar{u} \\
N_{2}(\bar{u})& = & \bar{u} 
\end{eqnarray}
For trimming curve I we compute the points along $\bar{u}=0$ as
\begin{eqnarray}
u_{I} (\bar{v})& = & \sum_{b=1}^{B} N_{b,p}^{I} (\bar{v})u_{b,I}\\
v_{I} (\bar{v})& = & \sum_{b=1}^{B} N_{b,p}^{I} (\bar{v})v_{b,I}
\end{eqnarray}
where $N_{b,p}^{I} (\bar{v})$ are the B-spline functions defining the trimming curve and $u_{b,I}$, $v_{b,I}$ are the local coordinates of the control points.
For triming curve II we compute the points at $\bar{u}=1$ as
\begin{equation}
u_{II} (\bar{v}) =  \sum_{b=1}^{B} N_{b,p}^{II} (\bar{v})u_{b,II}
\end{equation}
\begin{equation}
v_{II} (\bar{v}) =  \sum_{b=1}^{B} N_{b,p}^{II} (\bar{v})v_{b,II}
\end{equation}
where $N_{b,p}^{II} (\bar{v})$ are the B-spline functions defining the trimming curve and  $u_{b,II}$, $v_{b,II}$ are the local coordinates of the control points. This represents a linear interpolation between the trimming curves.
The evaluation of the integrals and the definition of the basis functions is carried out in the $\bar{u}$,$\bar{v}$ coordinate system and then mapped onto the $u,v$ and then the $x,y,z$ coordinate system (the mapping involves two Jacobians). An extension of the method to more than 2 trimming curves is possible.
The proposed mapping however would not work for the case where the trimming curve is a closed contour inside the $u,v$ domain. In this case the NURBS patch may be split into two or more patches.

\subsection{Geometry definition}

\begin{figure}
\begin{center}
\includegraphics[scale=1.0]{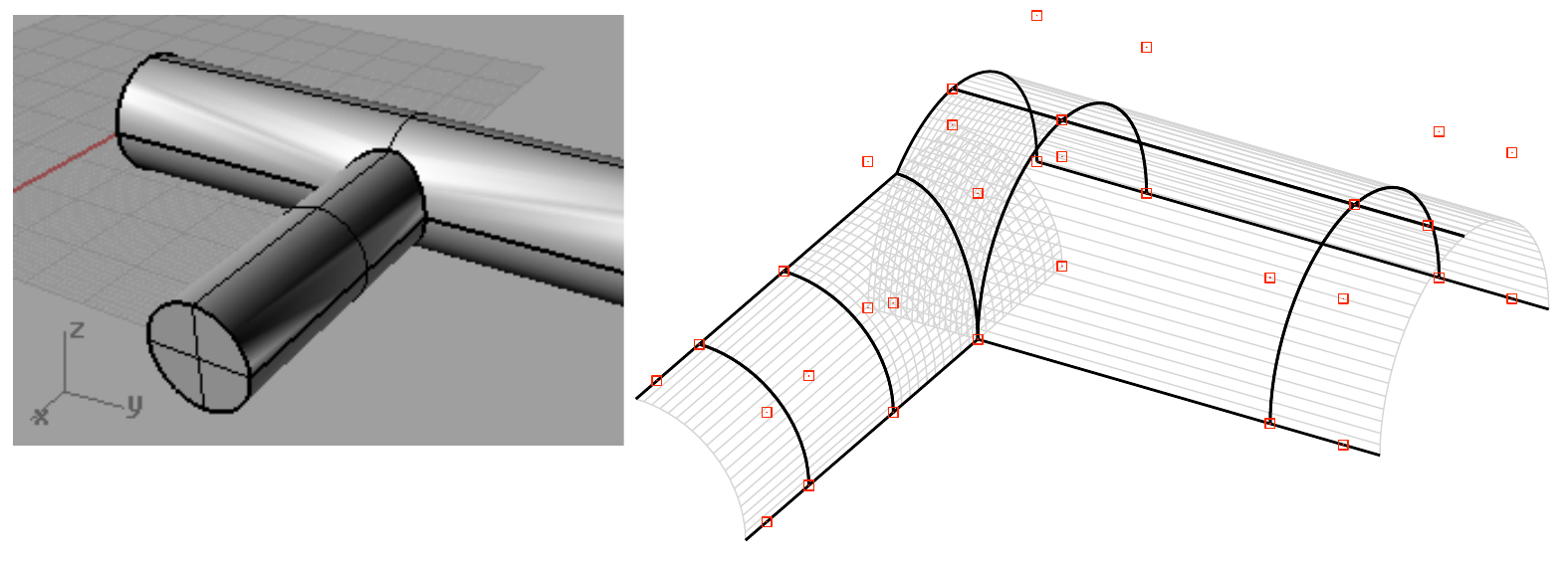}
\caption{CAD model of tunnel branch and resulting 1/4 simulation  model depicting NURBS patches and control points}
\label{Tun}
\end{center}
\end{figure}

To explain the definition of the problem geometry we use the example of a branched tunnel.
Figure \ref{Tun} depicts the CAD model of the tunnel with a branch at 90$^{\circ}$.
For the simulation, symmetry conditions were applied which meant that only 1/4 of the problem needed to be considered. The geometry definition with 6 NURBS patches (2 of them trimmed) and 3 infinite plane strain NURBS patches \cite{BeerBordas2014} is shown on the right of Figure \ref{Tun}. This geometry description is as accurate as the CAD description and needs no further refinement.

\section{SIMULATION }

An ideal companion to CAD is the boundary element method, as both rely on a description of the problem by surfaces. Therefore this method will be used for the simulation. In the following we only present a brief description of the method. Details of the implementation of the isogeometric BEM can be found in \cite{Scott2013197} and \cite{BeerBordas2014}.

\subsection{Boundary Element Method}
The boundary integral equation for an elastic continuum without body forces can be written as: 
\begin{eqnarray}
\label{eq:bie}
	\textbf{c}\left( P \right) \mathbf{u}\left( P \right) & = &
 \int_{S} \mathbf{U}\left( P,Q \right) \mathbf{t}\left( Q \right) dS  -   \int_{S}
\mathbf{T}\left( P,Q \right) \mathbf{u}\left( Q \right) dS 
\end{eqnarray}
The coefficient $ \mathbf{c}\left( P \right) $ is a free term
related to the boundary geometry. 
$ \mathbf{u}\left( Q \right) $ and $ \mathbf{t}\left( Q
\right) $ are the displacements and tractions on the boundary and
$\mathbf{U}\left( P,Q \right) $ and $ \mathbf{T}\left( P,Q \right) $  are matrices containing Kelvin's fundamental solutions (\textit{Kernels}) for the displacements and
tractions respectively.  $P$ is the source point and $Q$ is the field point. 

For the purpose of explaining the proposed simulation procedure we use the tunnel problem where tractions due to excavation are known and the displacements are unknown.
The integral equation can be discretized by using an interpolation of the displacements:
 \begin{equation}
\mathbf{u}=  \sum_{b=1}^{B}\sum_{a=1} ^{A}R_{a,b}^{p_{d},q_{d}} (\bar{u},\bar{v})\mathbf{d}_{a,b}^{e}
 \end{equation}
where  $\mathbf{d}_{a,b}^{e}$ denote the parameters for $\mathbf{u}$ at points $a,b$. The subscript $d$ of $p$ and $q$ indicates that the functions differ form the ones used for the description of the geometry. Therefore we use the terminology \textit{generalized IGA} as in the classical IGA reported in the literature the same functions are used.

To solve the discretized integral equation we use the \textit{Collocation} method, that is we satisfy it only at discrete points on the boundary $P_{n}$
The discretized integral equation can be written as:
\begin{equation}
\label{eq:bie_a}
\begin{split}
\mathbf{c}\left( P_{n} \right) \: \sum_{b=1}^{B}\sum_{a=1} ^{A}R_{a,b}^{p_{d},q_{d}} (\bar{u},\bar{v})\mathbf{d}_{a,b}^{ec}= 
\int_{S_{e}}
\mathbf{U}\left( P_{n},Q \right) \: \textbf{t }\: dS \\ - \sum_{e=1}^{E}
 \sum_{b=1}^{B}\sum_{a=1} ^{A}\int_{S_{e}} (\mathbf{T}\left( P_{n},Q \right) \: R_{a,b}^{p_{d},q_{d}} (\bar{u},\bar{v}) dS) \mathbf{d}_{a,b}^{e}\\ \quad for \quad n=1,2,3...N
\end{split}
 \end{equation}
where $ec$ denotes the patch that contains the collocation point and $E$ is the number of patches. The integrals over NURBS patches are computed using Gauss Quadrature. In some cases the patches have to divided into integration regions. For the case where the integrand tends to infinity inside a patch special procedures have to be applied. The reader is referred to \cite{BeerBordas2014} for a detailed discussion on this topic.

\begin{figure}
\begin{center}
\includegraphics[scale=0.6]{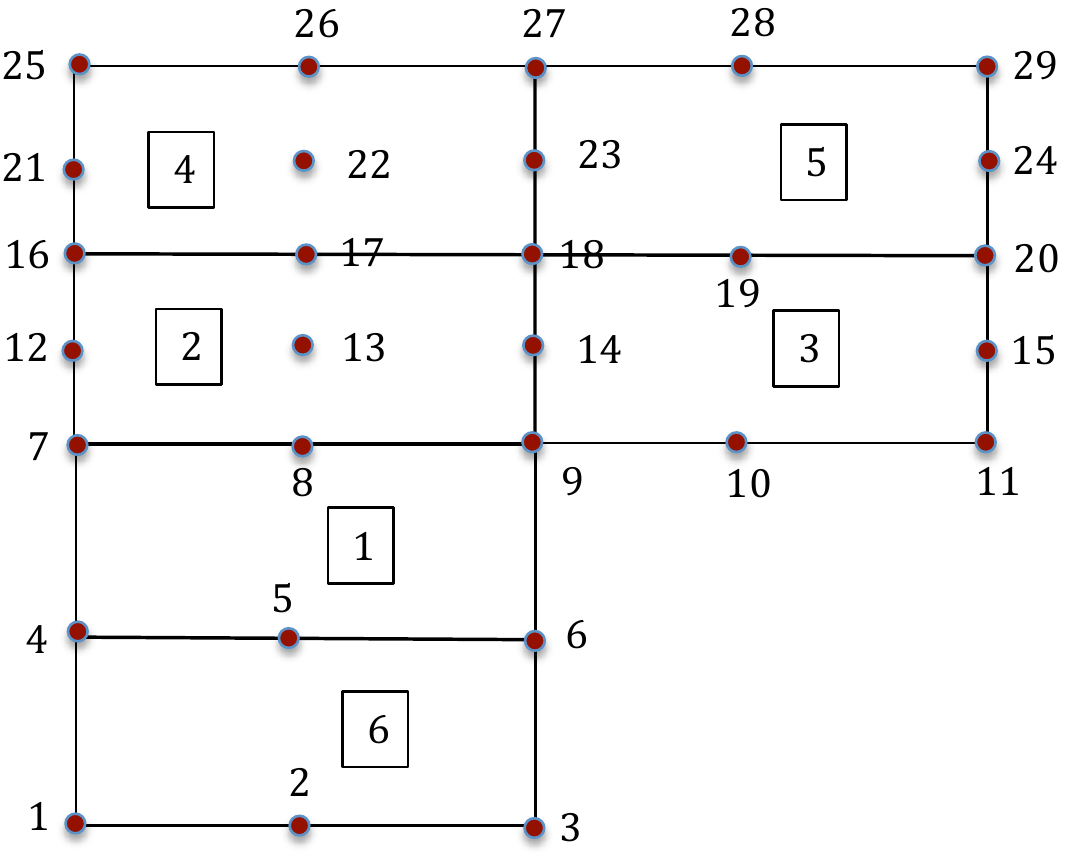}
\caption{Plot showing numbering of patches and parameters for tunnel example}
\label{DOF}
\end{center}
\end{figure}

The collocation points $P_{n} $ are first computed in the local coordinate system and then transferred to the global coordinate system as explained previously. 
The local coordinate of the collocation points can be computed using the method proposed by Greville \cite{Greville}:
\begin{equation}
\bar{u}(P_{i})= \frac{\bar{u}_{i+1}+\bar{u}_{i+2} + \dots +\bar{u}_{i+p_{d}}}{p_{d}} \qquad i=0,1, \dots ,I
\end{equation}
\begin{equation}
\bar{v}(P_{j})= \frac{\bar{v}_{j+1}+\bar{v}_{j+2} + \dots +\bar{v}_{j+q_{d}}}{q_{d}} \qquad j=0,1, \dots ,J
\end{equation}
where $i$ and $j$ denote the local (patch) numbering of the collocation point and $I$ and $J$ are the number of parameter points of the patch in $\bar{u}$ and $\bar{v}$ direction. $\bar{u}_{n}$ and $\bar{v}_{n}$ denote the corresponding entries in the Knot vector of the basis functions approximating the unknown displacements.

\textbf{Remark}: The Greville formulae also compute collocation points at the edges of NURBS patches. The coordinates of these points, computed in local coordinates of the different connecting NURBS, have to be the same. For trimmed NURBS this would only be the case if the parameter spaces of the trimming curves match. For the tunnel example presented here this was the case, so the collocation points matched, but this may not be guaranteed for a general application. However, this can be resolved by using discontiuous collocation \cite{Marus}.
The final system of equations to be solved is
\begin{equation}
\label{ }
[\textbf{T}] \{ \textbf{u} \}= \{ \textbf{F} \}
\end{equation}
where $[\textbf{T}]$ and $\{ \textbf{F} \}$ are assembled from patch contributions and $\{ \textbf{u} \}$ contains all displacement values.

\subsection{Solution and refinement strategies}

\begin{figure}
\begin{center}
\includegraphics[scale=0.8]{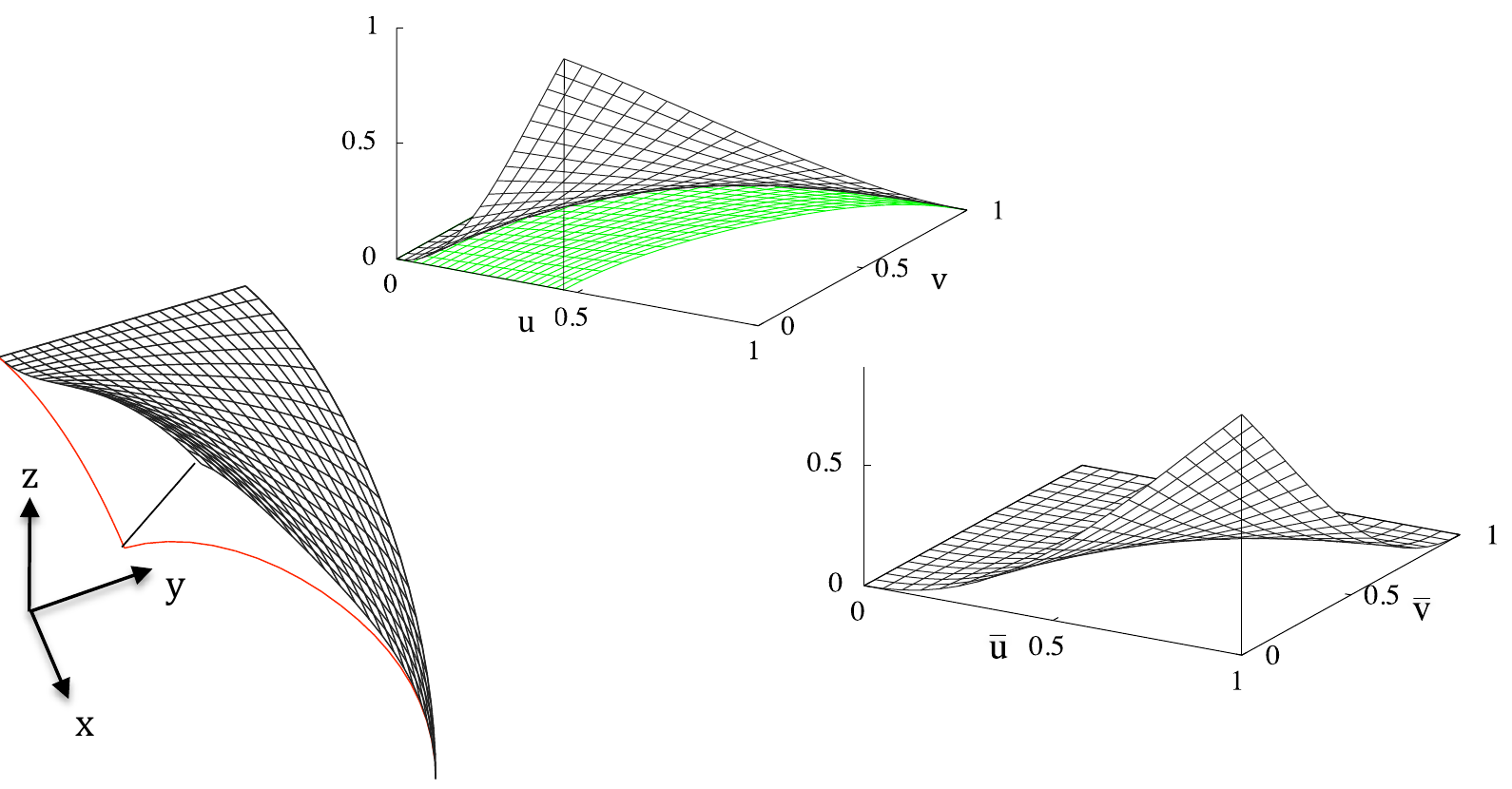}
\caption{Example of mapping a basis function onto patch 2}
\label{Mapf}
\end{center}
\end{figure}

\begin{figure}
\begin{center}
\includegraphics[scale=0.7]{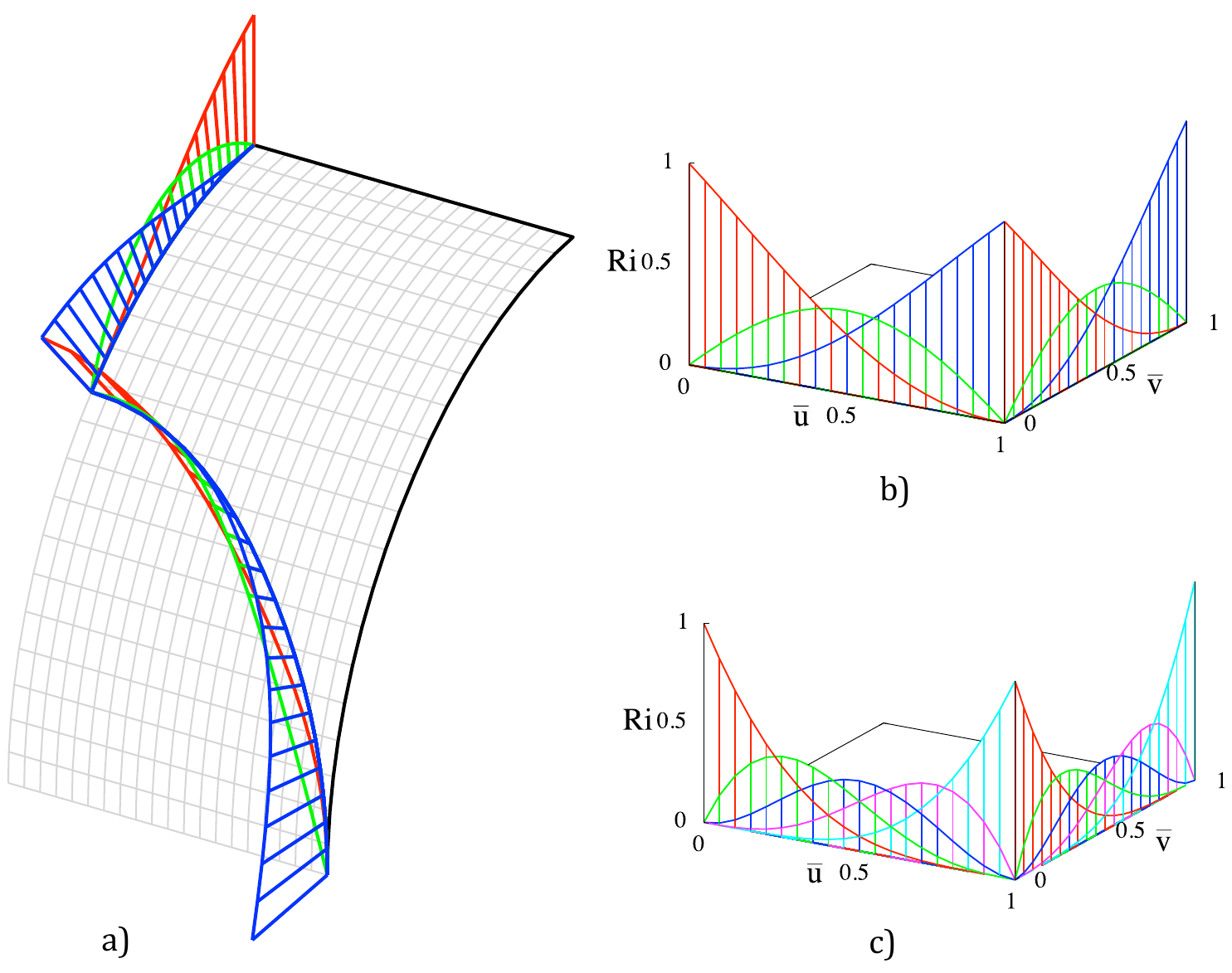}
\caption{Refinement process: Basis functions on edges of NURBS patch 2: a) in global coordinate system for $p=q=2$, b) in local coordinate system  for $p=q=2$ and c) refinement to $p=q=4$}
\label{Ref}
\end{center}
\end{figure}

Because the idea is to completely separate the description of the geometry from the approximation of the unknown  the unknown parameters are numbered independently. The parameter numbering is shown for the tunnel example in Figure \ref{DOF} for the coarsest discretization. This is updated automatically by the program during the refinement process.
The variation of the unknown is defined in the local $\bar{u},\bar{v}$ coordinate system and then mapped into the global coordinate system. For trimmed patches the procedure outlined in 2.2 is used. Figure \ref{Mapf} shows this for patch 2 of the tunnel example and the first basis function. We start the simulation with basis functions of order $p=q=2$ and use the different refinement strategies available in IGA. For the problem of the tunnel intersection it was found that raising the order to $p=q=4$ gave the best results (Figure \ref{Ref}).

\section{Example}

For testing the proposed algorithm we analyze the tunnel intersection with the following properties 
\begin{itemize}
  \item Elastic domain with E= 1000 MPa, v= 0
  \item Virgin stress: $\sigma_{z}=1 Mpa$ compression, all other components zero
  \item Symmetry about x-z and x-y planes
  \item Single stage excavation
\end{itemize}
For the case of a single stage excavation this is a pure Neumann problem and the tractions are given as:
\begin{equation}
\label{ }
\textbf{t}= \textbf{n} \boldsymbol \sigma
\end{equation}
where $\textbf{n}$ is the outward normal and $\boldsymbol \sigma$ is the pseudo-vector of virgin stress.

\begin{figure}
\begin{center}
\includegraphics[scale=0.4]{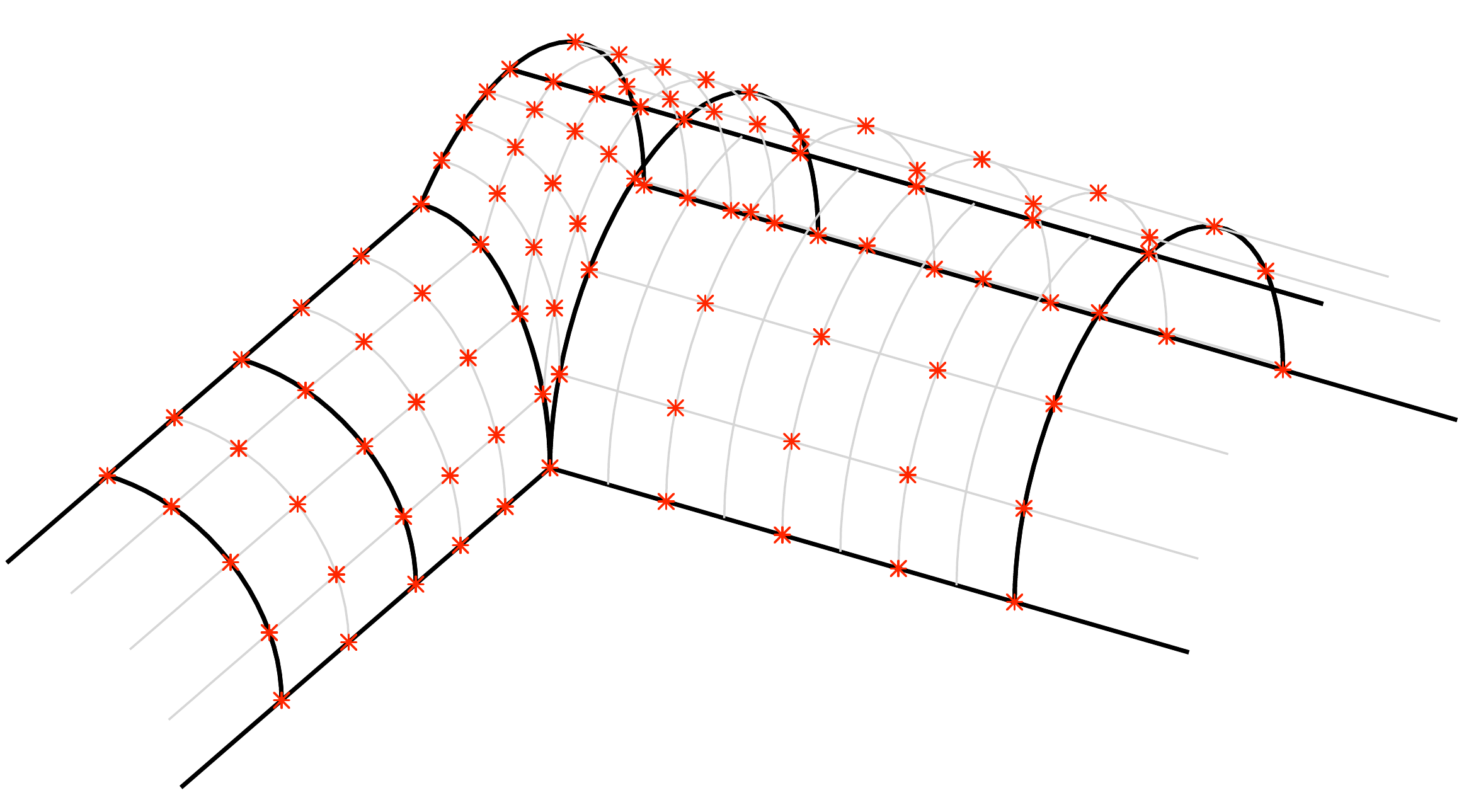}
\caption{Location of collocation points and subdivision into integration regions}
\label{Mesh}
\end{center}
\end{figure}

Figure \ref{Mesh} shows the location of collocation points and the subdivision into integration regions for the finest discretization ($p=q=4$, 291 degrees of freedom).  In the program, subdivision lines are generated automatically through collocation points. In addition, the user may define additional subdivision lines. Further subdivisions are automatically made by the program for the case where the source point $P$ is close to the integration region, using a quad tree method (for details see \cite{Beer2015})

\begin{figure}
\begin{center}
\includegraphics[scale=0.4]{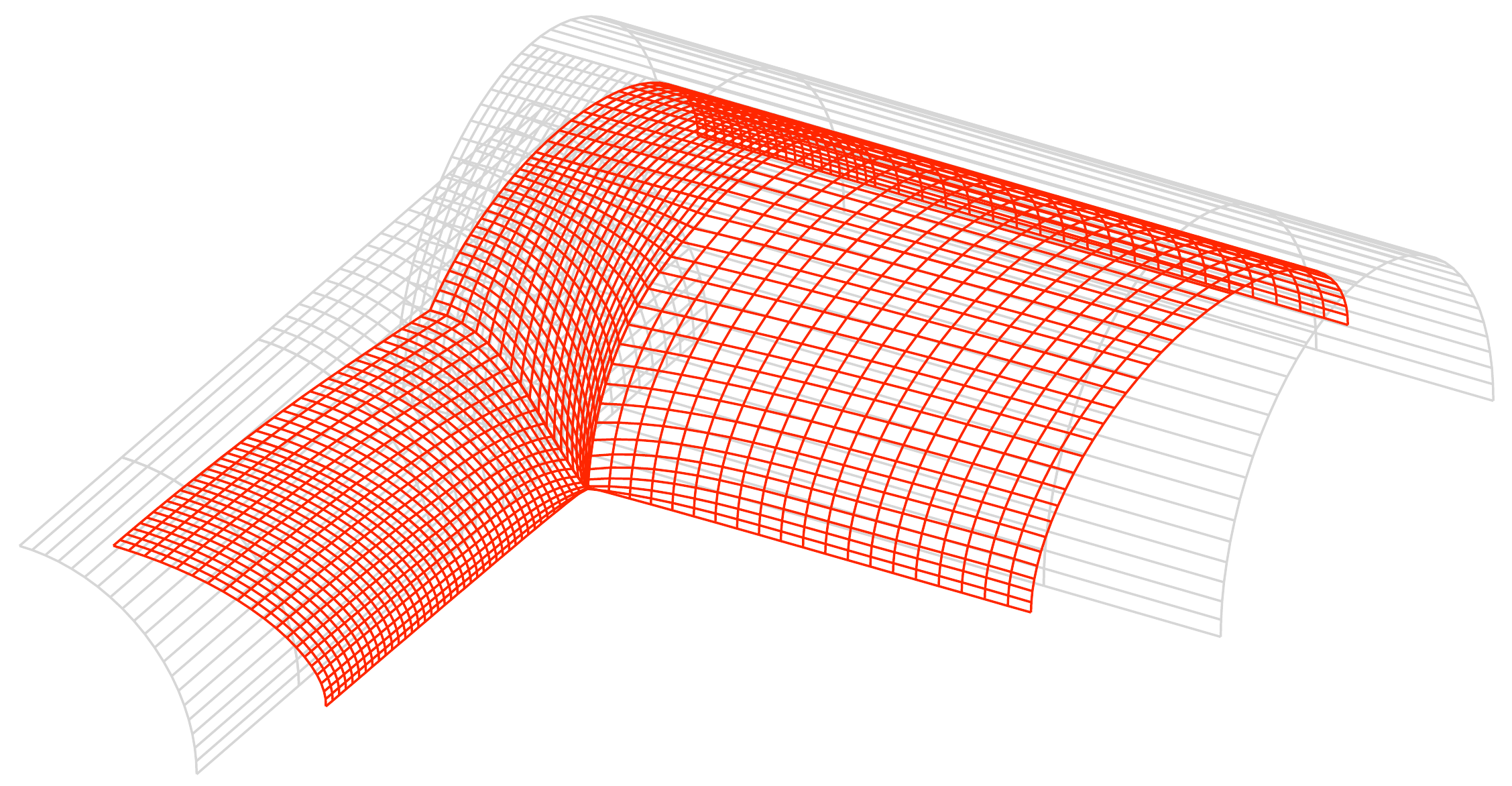}
\caption{Deformed shape}
\label{Def}
\end{center}
\end{figure}

\begin{figure}
\begin{center}
\includegraphics[scale=0.4]{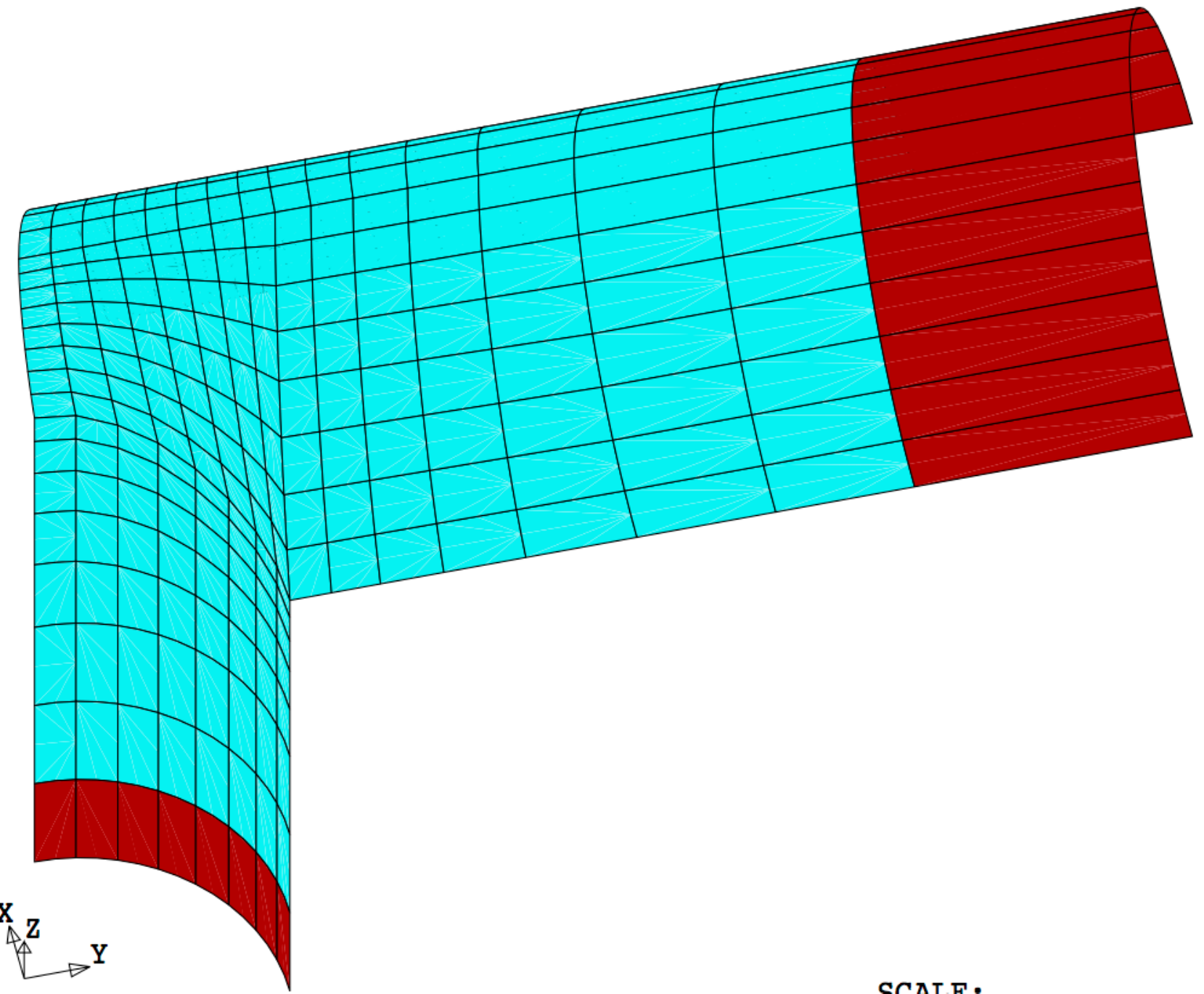}
\caption{Mesh used for the conventional BEM analysis with isoparametric elements}
\label{BEFE}
\end{center}
\end{figure}

\begin{figure}
\begin{center}
\includegraphics[scale=0.35]{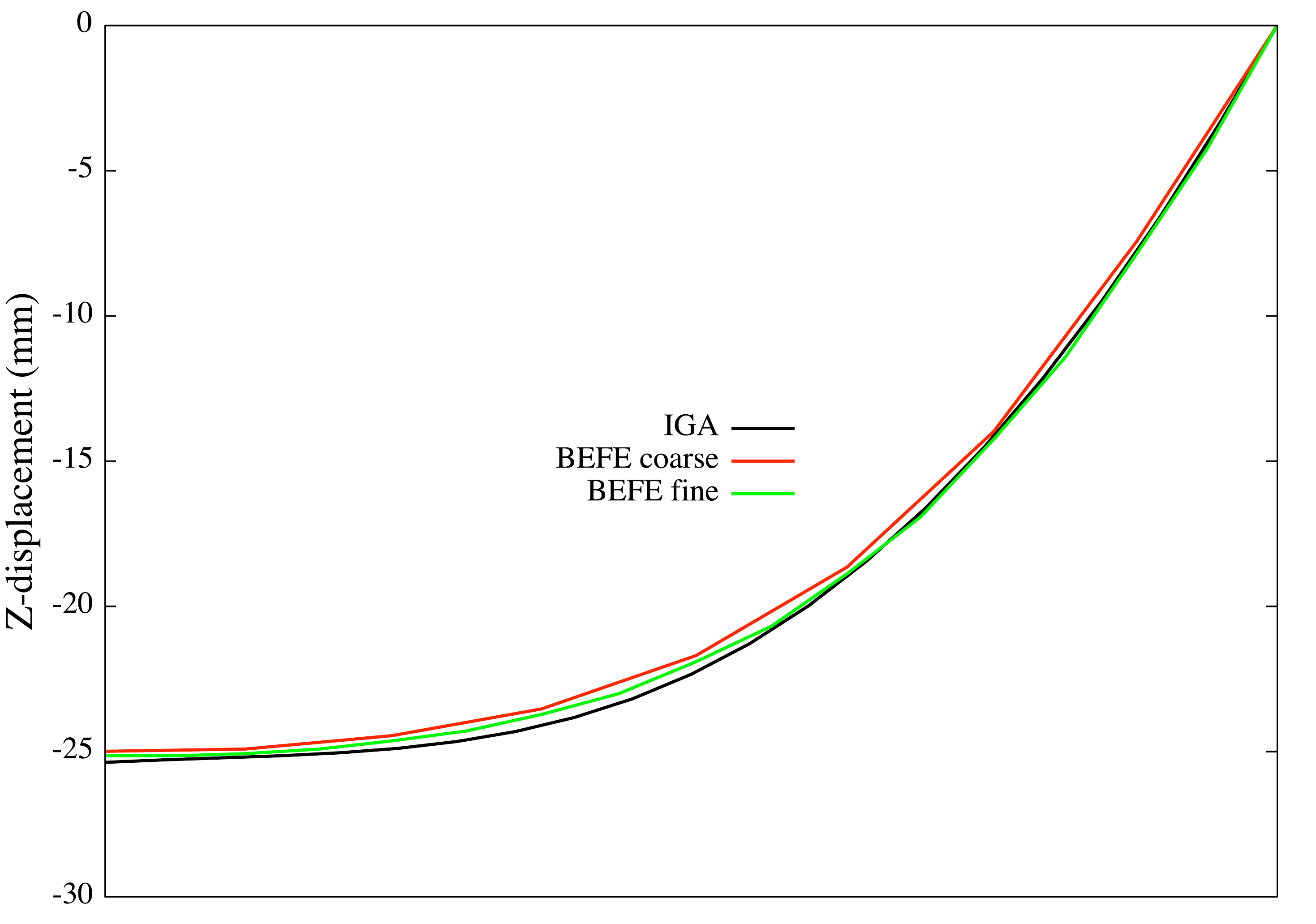}
\caption{Variation of the vertical displacement along the trimming line, comparison of new method (IGA) with isoparametric BEM (BEFE)}
\label{Comp}
\end{center}
\end{figure}

Figure \ref{Def} shows one result of the analysis namely the deformed shape.
To check the accuracy, the results are compared with a conventional BEM analysis using Serendipity functions for describing the geometry and the variation of the unknowns.
Figure \ref {BEFE} shows the mesh used for the analysis with the simulation program BEFE \cite{BEFE2003}. Two analyses were performed, one with linear and one with quadratic shape functions. The latter had 2895 unknowns.
The z-displacement along the trimming curve is shown in Figure \ref{Comp} for the conventional BEM and the new approach.
It can be seen that the conventional BEM results converge towards the results obtained by the method presented here.

\section{SUMMARY AND CONCLUSIONS}

We have presented a novel approach to the simulation with the boundary element method and trimmed NURBS patches. The innovations are in two parts. First, a procedure is presented for analyzing trimmed surfaces, which is much simpler to implement and more efficient than published methods. Secondly, we propose that the approximation of the unknown is completely independent from the definition of the boundary geometry. Our motivation comes from the fact that the boundary geometry is described with the same accuracy as the CAD model and needs no further refinement. Efficient refinement strategies available for NURBS  can then be applied to the description of the unknown only.

Comparison of the results of the analysis of a branched tunnel with a conventional isoparametric analysis shows good agreement. However, the number of degrees of freedom required to achieve the same (if not better) result is an order of magnitude smaller. The reason for this is that the approximation of the geometry is much more accurate and that the functions describing the variation of the unknowns exhibit a much higher continuity. In addition the need for the generation a mesh is completely avoided as data are taken directly from the CAD program. 

The algorithm for trimming still requires further testing on more complex cases and may have some limitations.
Future work will concentrate on the implementation of efficient procedures for non-linear analysis into the proposed framework and on methods of reducing storage and run times for large simulations.

\section*{ACKNOWLEDGMENTS}
This work was supported by the Austrian Science Fund (Project "Fast
isogeometric BEM"), Grant Number P24974-N30.

\bibliographystyle{plain} 
\bibliography{ifbbib}

\end{document}